\documentclass[11pt, a4paper]{article}

\usepackage{float}

\usepackage{times}

\usepackage{graphicx}

\usepackage{url}
\usepackage{amsmath}
\usepackage{amssymb}
\usepackage{cite}

\usepackage{booktabs}
\usepackage[table]{xcolor}

\usepackage{tikz}
\usetikzlibrary{arrows, petri, topaths}
\usepackage{tkz-berge}

\title{A Markovian model for association football possession and its outcomes}
\author{Javier L\'opez Pe\~{n}a \footnote{Kickdex Ltd, and Department of Mathematics, University College London. \texttt{javier@kickdex.com}}}
\date{}

\begin{document} \thispagestyle{empty}

\twocolumn
\maketitle

\thispagestyle{empty}
\begin{abstract}
	We propose a bottom-up approach to the study of possession and its outcomes for association football, based on probabilistic finite state automata with transition probabilities described by a Markov process. We show how even a very simple model yields faithful approximations to the distribution of passing sequences and chances of taking shots for English Premier League teams, which we fit using a whole season of granular game data (380 games). We compare the resulting model with classical top-down distributions traditionally used to describe possessions, showing that the Markov models yield a more accurate asymptotic behavior.
\end{abstract}

\section{Introduction}

In association football (football in the forthcoming) key game events such as shots and goals are very rare, in stark contrast to other team sports. By comparison, passes are two orders of magnitude more abundant than goals. It stands to reason that in order to get a comprehensive statistical summary of a football game passes should be one of the main focal points. However, not much attention has been paid to passing distributions, and the media attention to passes is normally limited to the total number of them, sometimes together with the \emph{passing accuracy}. A similar lack of attention is paid to the analysis of possession, oftentimes limited to a single percentage value.

In \cite{Reep1968} Reep and Benjamin analyze the distribution of length of passing sequences in association football and compare it to Poisson and negative binomial distributions. More recent works suggest the distribution of lengths of passing sequences in modern football is better explained by Bendford's law instead. However, all these works take a top-down approach in which a model is arbitrarily chosen and then fitted to the distribution, without any explanation on why football games should be described by that model.

In this work, we propose to take a bottom-up approach instead, inspired by our previous work on passing networks (cf. \cite{Lopez2013}) we model a team's game by a finite state automaton with states \textbf{Recovery}, \textbf{Possession}, \textbf{Ball lost}, \textbf{Shot taken}, with evolution matrix obtained from historical game data. The resulting Markovian system provides a model for possession after any given number of steps, as well as estimates for the probabilities of likelihood of either keeping possession or for it resulting in either of the two possible outcomes: lost ball or shot taken.

By choosing adequate fitting data, we compare the possession models for different teams and show how they vary across different leagues. We then compare the resulting model with the ones previously used in the literature (Poisson, NBD, Pareto) and with the actual possessions data in order to find the best fit. The obvious advantage of our model is that it explains the resulting distribution as the natural limit of a Markov process.

\section{Determining possession lengths}

Since football is such a fluid game, it is not always immediately clear which sequences of events constitute a possession. For our analysis, we consider the simplest possible notion of possession as any sequence of consecutive game events in which the ball stays in play and under control of the same team. As such, we will consider that a possession starts the first time a team makes a deliberate action on the ball, and ends any time the ball gets out of the pitch or there is a foul (regardless who gets to put the ball in play afterwards), any time there is a deliberate on-ball action by the opposing team, such as a pass interception tackle or a clearance (but not counting unsuccessful ball touches which don't interrupt the game flow), or any time the team takes a shot, regardless the outcome is a goal, out, hitting the woodwork, or a goalkeeper save.

One could introduce a further level of sophistication by distinguishing between clear passes and `\textit{divided balls}', such as passes to a general area where player of both teams dispute the ball (for instance in an aerial duel). For the sake of simplicity, we will not take this into consideration.

It is also worth noting that when measuring possession length we shall consider all in-ball actions, and not just passes. In particular dribbles, won aerial duels, and self-passes will all be considered as valid actions when accounting for a single possession.

\begin{figure}[ht]
    \includegraphics[scale=0.28]{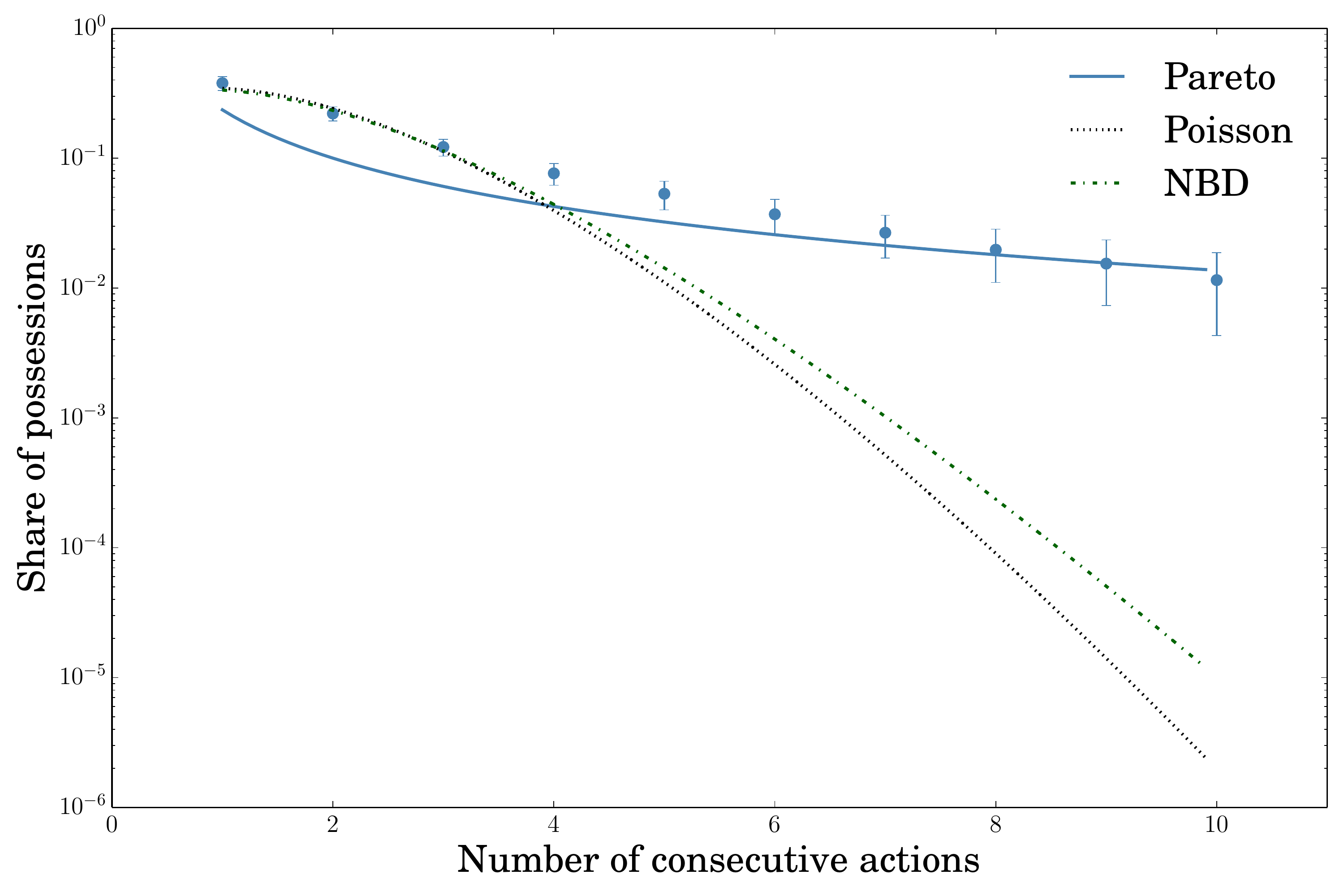}
    \caption{Classical top-down models}\label{fig:EPL0}
\end{figure}

In \cite{Reep1968} it is suggested that length of passing sequences can be approximated by Poisson or Negative Binomial distributions. Our tests suggest this is no longer the case for the average EPL team in the 2012/2013 season: Figure \ref{fig:EPL0} shows how the best fitting Poisson and NBD grossly underestimate the share of long possessions. For the sake of comparison, we have also included the fitting of a Pareto distribution, which displays the opposite effect (underestimating short possessions and overestimating longer ones). For our data, all three classic distributions fail to meet the asymptotic behavior of the observed data.

This can be partly explained by our different notion for possession length, since Reep and Benjamin only consider the total number of passes in a possession, but the different definition does not tell us the whole story. According to Reep and Benjamin data, there are only 17 instances, measured over 54 games in 1957-58 and 1961-61, of sequences involving 9 passes or more (out of around 30000 possessions). Even if we use their more restrictive notion of possession there are many more instances of long possession nowadays. A possible explanation might be the generally admitted fact that football playing style has evolved over the years, leaning towards a more and more technical style which favors longer possessions, with many teams consistently playing possessions of 20 passes or even more.

In any case, it is quite evident that in order to accurately describe possessions in current professional football we need to find a different type of model which allows for a longer tail. One possible such candidate would be a power-law distribution, such as Pareto's (also plotted in Figure \ref{fig:EPL0} for comparison). As we plot longer parts of the tail, however, it will become apparent that Pareto distribution does not display the correct type of asymptotic behavior in order to describe our data.

\section{Markovian possession models}

A typical possession in a football game starts with a ball recovery, which can be achieved in an active or passive manner. Ball recovery is followed by a sequence of ball movements that will eventually conclude with a loss of possession. Said lost of possession can be either intentional, due to an attempt to score a goal, or unintentional, due to an error, an interception/tackle, or an infringement of the rules.

The different possession stages can be modeled using a nondeterministic finite state automaton, with initial state indicating the start of a possession, final states indicating the end of the possession, and some number of other states to account for intermediate stages of the possession, as well as appropriate probabilities for the state transitions.

The freedom in the choice of intermediate states allows for a wealth of different Markov models with varying degrees of complexity. One of the simplest such possible schemes is to consider only one intermediate state ``\textbf{Keep}'' to indicate continued possession, and two different final states ``\textbf{Loss}'' and ``\textbf{Shot}'' indicating whether the possession ends in a voluntary manner (by taking an attempt to score a goal) or in an involuntary manner. This simple model is schematized by the following diagram:

\begin{figure}[ht]
    \begin{tikzpicture}[scale=0.8,transform shape, >=stealth']
        \Vertex[x=0, y=0]{Recovery}
        \Vertex[x=4, y=0]{Keep}
        \Vertex[x=8, y=1]{Loss}
        \Vertex[x=8,y=-1]{Shot}
        \tikzset{EdgeStyle/.style={->}}
        \Edge[style= dashed](Recovery)(Keep)
        \Edge[label = $p_k$, style = loop](Keep)(Keep)
        \Edge[label = $p_l$](Keep)(Loss)
        \Edge[label = $p_s$](Keep)(Shot)
    \end{tikzpicture}
    \caption{Simple Markov model}
\end{figure}
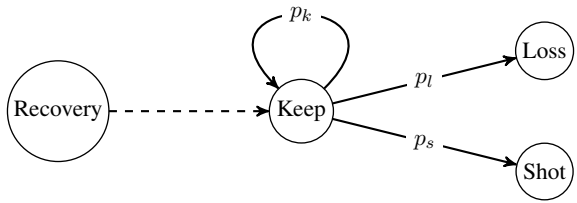

A slight variation allows to track for divided balls (though using this model would require us to change the definition of possession that we described above):

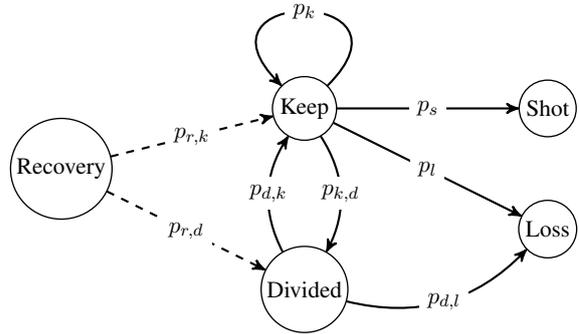
\begin{figure}[ht]
    \begin{tikzpicture}[scale=0.8,transform shape, >=stealth']
        \Vertex[x=0, y=0]{Recovery}
        \Vertex[x=4, y=1]{Keep}
        \Vertex[x=8, y=-1]{Loss}
        \Vertex[x=8,y=1]{Shot}
        \Vertex[x=4, y=-2]{Divided}
        \tikzset{EdgeStyle/.style={->}}
        \Edge[style=dashed, label=$p_{r,k}$](Recovery)(Keep)
        \Edge[style=dashed, label=$p_{r,d}$](Recovery)(Divided)
        \Edge[label=$p_k$, style=loop](Keep)(Keep)
        \Edge[label=$p_l$](Keep)(Loss)
        \Edge[label=$p_s$](Keep)(Shot)
        \Edge[label=$p_{k,d}$, style=bend left](Keep)(Divided)
        \Edge[label=$p_{d,k}$, style=bend left](Divided)(Keep)
        \Edge[label=$p_{d,l}$, style=bend right](Divided)(Loss)
    \end{tikzpicture}
    \caption{Markov model with divided balls}
\end{figure}

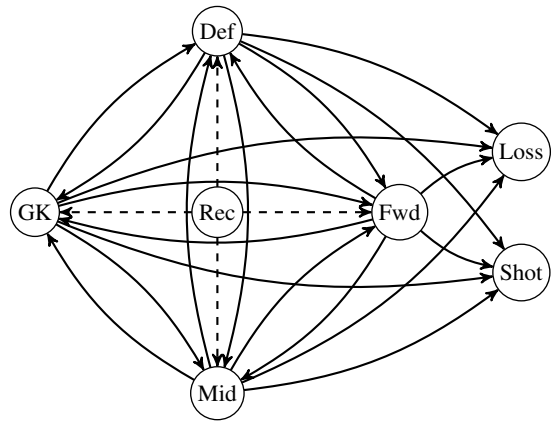
\begin{figure}[ht]
    \begin{tikzpicture}[scale=0.8,transform shape, >=stealth']
        \Vertex[x=3, y=0]{Rec}
        \Vertex[x=0, y=0]{GK}
        \Vertex[x=3, y=3]{Def}
        \Vertex[x=3, y=-3]{Mid}
        \Vertex[x=6, y=0]{Fwd}
        \Vertex[x=8, y=1]{Loss}
        \Vertex[x=8,y=-1]{Shot}

        \tikzset{EdgeStyle/.style={->}}
        \Edge[style=dashed](Rec)(GK)
        \Edge[style=dashed](Rec)(Mid)
        \Edge[style=dashed](Rec)(Def)
        \Edge[style=dashed](Rec)(Fwd)
        \Edge[style={bend left = 15}](GK)(Def)
        \Edge[style={bend left = 15}](GK)(Mid)
        \Edge[style={bend left = 15}](GK)(Fwd)
        \Edge[style={bend left = 15}](Def)(GK)
        \Edge[style={bend left = 15}](Def)(Mid)
        \Edge[style={bend left = 15}](Def)(Fwd)
        \Edge[style={bend left = 15}](Mid)(GK)
        \Edge[style={bend left = 15}](Mid)(Def)
        \Edge[style={bend left = 15}](Mid)(Fwd)
        \Edge[style={bend left = 15}](Fwd)(GK)
        \Edge[style={bend left = 15}](Fwd)(Mid)
        \Edge[style={bend left = 15}](Fwd)(Def)

        \Edge[style={bend left = 15}](GK)(Loss)
        \Edge[style={bend right = 15}](GK)(Shot)
        \Edge[style={bend left = 15}](Def)(Loss)
        \Edge[style={bend left = 15}](Def)(Shot)
        \Edge[style={bend right = 15}](Mid)(Loss)
        \Edge[style={bend right = 15}](Mid)(Shot)
        \Edge[style={bend left = 15}](Fwd)(Loss)
        \Edge[style={bend right = 15}](Fwd)(Shot)
    \end{tikzpicture}
    \caption{Positional Markov model}
    \label{dia:positional}
\end{figure}

In an alternative model, one might be interested on the interactions among players in different positions (loops have been removed from the graph to avoid cluttering). This model is sketched in Figure \ref{dia:positional}.

Once we have chosen the set of states that configure our Markov model, its behavior is then described by the \emph{transition matrix} $A=(a_{i,j})$, where entry $a_{i,j}$ represents the probability of transitioning from state $i$ to state $j$; similarly, given an initial state $i$, the probability of the game being at state $j$ after $r$ steps is given by $A^r e_i$, where $e_i = (0, \dotsc, 0, 1, 0, \dotsc, 0)$ is the vector with coordinate 1 in position $i$ and $0$ everywhere else.

More complex Markov models allow for a richer representation of game states, but since transition probabilities must (in general) be determined heuristically, they will be harder to fine-tune.

\section{Model fitting}

As a proof of concept, we will perform a detailed fitting of the simplest Markov model described in the previous section. This first approximation model works under the very simplistic assumption that the transition probabilities between two game states remain constant through the entire sequence of events. Transition probabilities $p_k, p_l, p_s$ (which must satisfy $p_k + p_s + p_l = 1$) are, respectively, the probability of keeping possession, losing the ball, or taking a shot.

Since the transition matrix for this system is extremely simple, this model can be solved analytically, yielding the probability distribution given by
\begin{equation}
    P(\{X = x\}) = (1-p_k)p_k^{x-1} \simeq Ce^{-\lambda x},
\end{equation}
(where $\lambda = -\log p_k$ and $C$ is a normalization constant) suggesting that the distribution of lengths follows a pattern asymptotically equivalent to an exponential distribution.

\begin{figure}[ht]
    \includegraphics[scale=0.28]{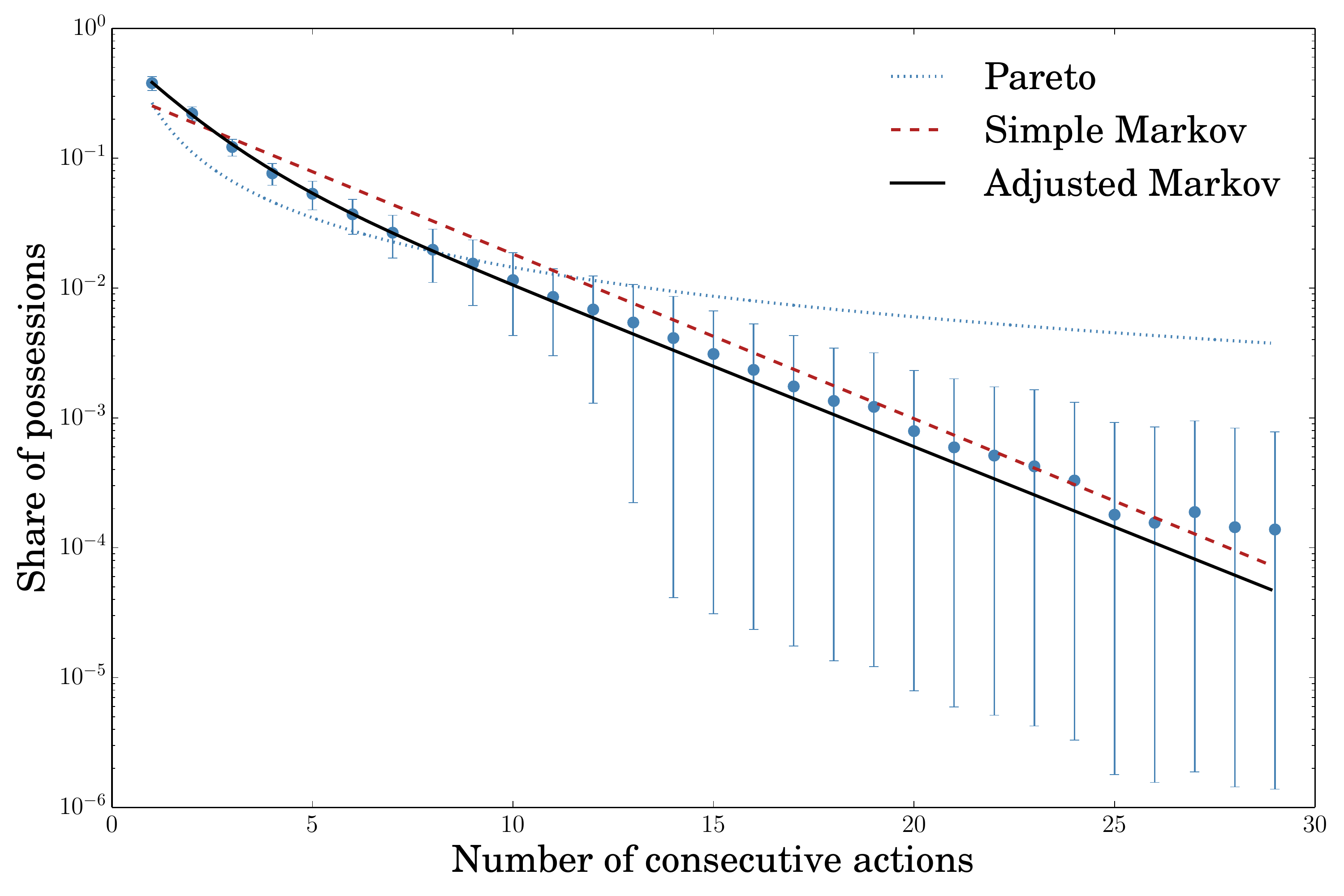}
    \caption{Power law and Markov models}\label{fig:EPL1}
\end{figure}

As Figure \ref{fig:EPL1} shows, the simple Markov model provides a good fitting model for the general trend, with the correct asymptotic behavior, but it tends to over-estimate occurrence of sequences of 3 to 8 actions, as well as severely underestimating the number of sequences consisting on a single action. For individual Premier League teams, Pearson's goodness of fit test yields a $p$--value higher than $0.99$ in all cases.

The simple Markov model can be easily improved by weakening the constraint on the transition probability being constant. A simple way of doing it is by adding an accumulation factor $b$ (with $0<b<1$) and modifying the probability distribution to
\begin{equation}
    P(\{X = x\}) = C(e^{-\lambda x} + b^x),
\end{equation}
since $b^x$ goes to 0 as $x$ increases, the added factor does not modify the asymptotic behavior of the resulting probability distribution, but it allows to correct for the errors in the probabilities for short possessions. It is worth noting that whilst this adjusted model (also shown in Figure \ref{fig:EPL1}) does not strictly come from a Markov process (as the transition probability is no longer constant), the additional factor $b$ can be interpreted as a (negative) self-affirmation feedback, incorporated to the model in a similar way as the one used in \cite{Bittner2007} for goal distributions. A possible interpretation of this factor would be the added difficulty of completing passes as a team proceeds to move the ball closer to their opponents box. A higher value of the $b$ coefficient means that a team is more likely to sustain their passing accuracy over the course of a long possession.

\begin{figure}[ht]
    \includegraphics[scale=0.28]{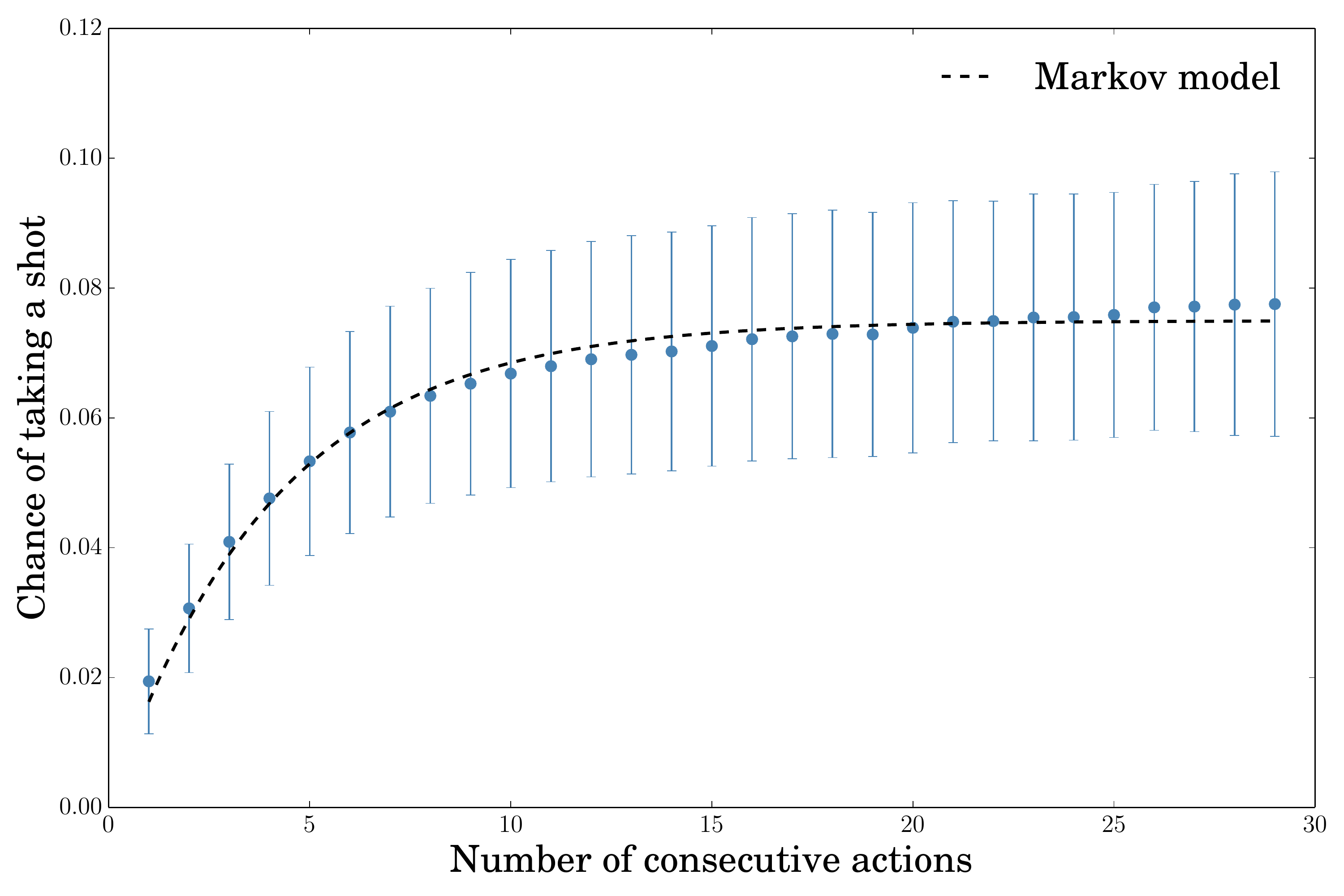}
    \caption{Probability of shots for the Markov model}\label{fig:shotchance}
\end{figure}

Besides the distribution of length of possessions, the Markov model allows us to study the probability of any state of the model after a given number of steps. In particular, one can look at the `\textbf{Shot}' state, obtaining a model for the probability that a shot will have taken place within a given number of actions. Once again, the system can be solved analytically, yielding
\begin{equation}
    P(\{\text{\textbf{Shot}} | X\leq x \}) = \sum_{i=0}^{x-1} p_k^{i}p_s = p_s \frac{1-p_k^x}{1-p_k},
\end{equation}
where this probability should be understood as the chance that a team will take a shot within $x$ consecutive ball touches. Figure \ref{fig:shotchance} shows once again how the model yields a good approximation to the observed data fitting comfortably within the error bars. The $p$--value for the Pearson's goodness of fitting test is again greater than $0.99$.

The fitted coefficients for all the  Premier League teams in the 2012-2013 season are listed in table \ref{tbl:markovfitting} (teams are sorted by final league position). There is a remarkable correlation between higher values of the `\emph{keep probability}' $p_k$ and what is generally considered `\emph{nice gameplay}', as well as between higher values of the shot probability $p_s$ and teams considered to have a more aggressive football style.

Traditionally possession is considered a bad indicator of performance, but our number suggest that this does not need to be the case when it is analyzed in a more sophisticated manner. Apart from a few outliers, such as under-performers Swansea and Wigan Athletic, and over-performing Stoke, most teams fitted probabilities correlate fairly well with their league table positions.

\begin{table}[t]
\begin{tabular}{rlccc}
\toprule
 & Team &  $p_k$ &  $p_s$ &  $b$ \\
\midrule
1   & Manchester United &  0.794 &  0.017 &  0.685 \\
2   & Manchester City &  0.794 &  0.016 &  0.683 \\
3   & Chelsea &  0.785 &  0.018 &  0.663 \\
4   & Arsenal &  \textbf{0.797} &  0.015 &  \textbf{0.690} \\
5   & Tottenham Hotspur &  0.782 &  0.016 &  0.653 \\
6   & Everton &  0.772 &  0.017 &  0.631 \\
7   & Liverpool &  0.788 &  \textbf{0.019} &  0.672 \\
8   & West Bromwich Albion &  0.771 &  0.018 &  0.630 \\
9   & Swansea City &  0.794 &  0.018 &  0.684 \\
10  & West Ham United &  0.756 &  0.018 &  0.589 \\
11  & Norwich City &  0.764 &  0.014 &  0.601 \\
12  & Fulham &  0.783 &  0.018 &  0.657 \\
13  & Stoke City &  0.752 &  0.015 &  0.575 \\
14  & Southampton &  0.772 &  0.015 &  0.630 \\
15  & Aston Villa &  0.771 &  0.016 &  0.619 \\
16  & Newcastle United &  0.771 &  \textbf{0.019} &  0.624 \\
17  & Sunderland &  0.763 &  0.017 &  0.601 \\
18  & Wigan Athletic &  0.783 &  0.017 &  0.660 \\
19  & Reading &  0.749 &  0.016 &  0.558 \\
20  & Queens Park Rangers &  0.767 &  0.016 &  0.618 \\
\bottomrule
\end{tabular}
\caption{EPL Markov model fitting parameters}
\label{tbl:markovfitting}
\end{table}

\section{Conclusions and future work}

We have shown how Markov processes provide a bottom-up approach to the problem of determining the probability distribution of possession related game states, as well as their outcomes. The bottom up approach has an obvious advantage of providing a conceptual explanation for the resulting probability distribution, but besides that we have shown that even a very simple Markov model yields better approximations than the classical top-down approaches.

It is worth noting that even if we focused in the particular case of association football (and concretely the English Premier League) the type of analysis we have carried out is of a very general nature and can easily be performed for many team sports which follow similar possession patterns, such as basketball, hockey, handball, or waterpolo, to name a few.

Besides the study of probability distributions, sufficiently granular Markov processes can be used to carry out game simulations. In theory one might want to consider a full Markov model containing a state for every single player, as in the passing networks described in \cite{Lopez2013}, and use it as the base of an agent based model in order to forecast game outcomes, although finding all the probability transitions in that extreme situation would admittedly be very hard.

\section{Data and analysis implementation details}

Our analysis uses data for all the English Premier League games in the 2012-2013 season (380 games). Raw data for game events was provided by Opta. Data munging, model fitting, analysis, and chart plotting were performed using IPython \cite{ipython} and the python scientific stack \cite{numpy, matplotlib}.

\end{document}